\documentstyle[12pt,leqno]{article}
\renewcommand{\P}{\mbox{${\bf P}$}}
\newcommand{\Schubert}{{\sc schubert}}
\renewcommand{\O}{\mbox{${\cal O}$}}
\newcommand{\Q}{\mbox{${\cal Q}$}}
\newcommand{\C}{\mbox{${\bf C}$}}
\renewcommand{\l}{\lambda}

\newcommand{\s}[1]{\mbox{$\sigma_{#1}$}}
\newtheorem{fact}{Theorem}

\begin{document}
\bibliographystyle{plain}

\begin{flushright}
OSU-M-92-3
\end{flushright}

\bigskip
\begin{center}
Rational curves on Calabi-Yau manifolds: verifying\\
predictions of
Mirror Symmetry\\
\bigskip
Sheldon Katz\\
Department of Mathematics\\
Oklahoma State University\\
Stillwater, OK 74078\\
email: katz@math.okstate.edu
\end{center}

\section*{}
Recently, mirror symmetry, a phenomenon in superstring theory,
has been used to give tentative calculations of several numbers
in algebraic geometry \nolinebreak\footnote{See the papers in \cite{yau} for
general background on mirror symmetry.}.  This yields predictions for the
number of rational
curves of any degree $d$ on general Calabi-Yau hypersurfaces in
$\P^4\ \cite{cogp},\ \P(2,1^4),\ \P(4,1^4)$, and $\P(5,2,1^3)$
\cite{fontper,kt,morp-f}.
The techniques used in the calculation rely on
manipulations of path integrals which have not yet been put on a rigorous
mathematical footing.  On the other hand,
there is currently
no prospect of calculating most of these numbers by algebraic geometry.

Until this point, three of these numbers have been verified, all for the
quintic hypersurface in $\P^4$: the number of lines (2875) was known
classically, the number of conics (609250) was calculated in
\cite{finite}, and the number of twisted cubics (317206375)
was found recently by
Ellingsrud and Str{\o}mme \cite{escub}.

Even more recently \cite{gmp}, higher dimensional mirror symmetry has been
used to
predict the number of rational curves on Calabi-Yau hypersurfaces in higher
dimensional projective spaces which meet 3 linear subspaces of certain
dimensions.  Again, there is no known way to calculate these using algebraic
geometry.

The purpose of this paper is to verify some of these numbers in low degree,
giving more evidence for the validity of mirror symmetry.  In \S 1, the
number of weighted lines in a weighted sextic in $\P(2,1^4)$ is calculated,
as well as the number of weighted lines in a weighted octic in $\P(4,1^4)$.
In \S 2, the number of lines on Calabi-Yau hypersurfaces of dimension up to
10 which satisfy certain incidence properties is calculated.  In \S 3,
the number of conics on these same Calabi-Yau hypersurfaces satisfying the
same incidence properties is calculated.  These numbers are closely related to
the Gromov-Witten invariants defined in \cite{wittsm,morht,gmp}; and it is {\em
these}
numbers that are recorded here.  In all instances, the calculations
agree with those predicted by mirror symmetry.  Thus the number of verified
predictions has increased from 3 to 65.

There are two parts to all of these calculations.  The first part is to
express the desired numbers in terms of the standard constructs of intersection
theory.  The second part is to evaluate the number
using the Maple package \Schubert\
\cite{schub} (although the number of weighted lines in a weighted octic in
$\P(4,1^4)$ was first found via classical enumerative geometry, using a
classical enumerative formula).  The short \Schubert\ code is not included
here, but is available upon request.

While it is checked that the data being enumerated is finite, no attempt has
been made here to check that the multiplicities are 1.  All enumeration
takes multiplicities into consideration.  This suffices for comparison to
the numbers arising in
physics, since the Feynman path integrals would take account of any
multiplicities greater than 1 as well.

Some of the Gromov-Witten invariants were computed in \cite{gmp}
using an intriguing relation between the various invariants.  These
relations arise in
conformal field theory.  A mathematical proof of the relations for the
invariants corresponding to lines is sketched here.

It is appropriate to point out the recent work of Libgober and Teitelbaum
\cite{ltcy}, who have apparently correctly guessed the mirror manifold of
complete intersection Calabi-Yau threefolds in an ordinary projective space.
Their conjectured mirrors yield predictions for the numbers of rational
curves.  The predicted number of lines coincides with the results of a
calculation done by
Libgober 20 years ago, and the predicted
number of conics coincides with the results of an unpublished calculation
done by Str\o mme and Van
Straten in 1990.

I'd like to thank D.~Morrison for helpful suggestions and conversations,
and for encouraging
me to write this paper.  I'd also like to thank S.~Kleiman for his
suggestions which have improved the manuscript.

\section{Weighted projective spaces and their Grassmannians}

Let $\P(k,1^n)$ denote an $n$-dimensional weighted projective space with
first coordinate having weight $k>1$, all other coordinates having weight 1.
Thus $\P(k,1^n)$ consists of all non-zero $(n+1)$-tuples $(x_0,\ldots,x_n)$,
with $(x_0,\ldots,x_n)$ identified with $(\l^kx_0,\l x_1,\ldots,\l x_n)$ for
any $\l\neq 0$.
Note
that $\P(k,1^n)$ is smooth outside the singular point $p=(1,0,\ldots,0)$.
There is a natural rational projection map $\pi:\P(k,1^n) ---> \P^{n-1}$
defined outside $p$ given by omitting
the first coordinate.
Let
$X$ be a weighted hypersurface of weight $d$.  Assume $p\not\in X$ (this
implies that $k|d$, and that the monomial $x_0^{d/k}$ occurs in an
equation for $X$).  It is further assumed that $X$ is smooth.  The general
weighted hypersurface whose weight is a multiple of $k$ is an example of such
an $X$.

\bigskip\noindent
{\bf Definition}
A {\em weighted $r$-plane in $\P(k,1^n)$} is the image of
a section of $\pi$ over an
$r$-plane in $\P^{n-1}$.

\bigskip
Note that weighted $r$-planes do not contain $p$.

Let $P$ be a weighted $r$-plane, with $L$ its image in $\P^{n-1}$.  Let
$(q_0,\ldots ,q_r)$ be any homogeneous coordinates on $L\simeq\P^r$.  Then
$P$ may be thought of as the image of $L$ via the mapping
$x_0=f_k(q_0,\ldots,q_r),\ x_i=l_i(q_0,\ldots,q_r)$, where $f_k$ is a form of
degree $k$, and the $l_i$ are all linear.  Once $L$ is fixed,
we may fix in mind a choice of the $q_i$ and $l_i$.

The moduli space of weighted $r$-planes can be represented (and compactified)
as follows.  Conventions have been chosen to be consistent with those in
\Schubert\ \cite{schub}.  Let $G=G(r+1,n)$ be the Grassmannian of
$r$-dimensional linear
subspaces $L$
of $\P^{n-1}$.  Let $V$ be the $n$-dimensional vector space of linear forms
on $\P^{n-1}$.
This identifies $\P^{n-1}$ with $\P(V)=Proj(S^*V)$.  $G$ is then the space of
$r+1$
dimensional
quotients of $V$ (since the space of linear forms restricted to $L$ is an
$r+1$-dimensional quotient of $V$).  Let $Q$ be the universal rank $r+1$
quotient bundle on $G$.

The equation $x_0-f_k(q_0,\ldots,q_r)=0$ which describes a section of
$\pi$ over
an $r$-plane may be identified with a section $s$ of the bundle
$\C\oplus S^k(Q)$,
where \C\ denotes the trivial bundle.  A scalar multiple of this section
would correspond to the equation $ax_0-af_k(q_0,\ldots,q_r)=0$, which
defines the same weighted $r$-plane.
Note that (up to scalar) $s$ does
not depend on any of the choices which have been made.
So $M=\P(\C\oplus S^k(Q)^*)$ gives a compactification of the space of weighted
$r$-planes.  Here, $\P(E)$ denotes the space of rank 1 quotients of the fibers
of the bundle $E$; hence the need for dualizing in defining $M$.

In the sequel, we will also refer to $M^o\subset M$, the open subset which
corresponds to the actual
weighted $\P^1$'s, in other words, $M^o=M-\P(S^kQ^*)$, where $\P(S^kQ^*)$ is
included in $M$ via the map induced by the natural projection
$\C\oplus S^kQ^*\to S^kQ^*$.

\section{Lines on weighted hypersurface Calabi-Yau threefolds}
\label{weight-section}

Now let $k>1$, and let $X\subset\P(k,1^4)$ be a smooth weighted
hypersurface of weight
$k+4$ with $(1,0,0,0,0)\not\in X$.  As has been noted in the previous section,
this implies that $k|k+4$, which in turn implies that $k=2$ or $k=4$.  The
weight $k+4$ has been chosen to ensure that $X$ is Calabi-Yau, i.e. that $X$
has trivial canonical bundle.

The rational projection map $\pi$ restricts to a morphism
$\pi:X\to\P^3$.  This is
a 3-1 cover for $k=2$, and a 2-1 cover for $k=4$.  The goal of this section
is to enumerate the weighted $\P^1$'s contained in $X$.

Let us first consider the case $k=4$.  Then an equation for $X$ has the
form
\begin{equation}
\label{hyper}
F=ax_0^2+g_4(x_1,\ldots,x_4)x_0+g_8(x_1,\ldots,x_4)=0,
\end{equation}
where $a\in\C$
and
$g_i$ has degree $i$ for $i=4$ or 8.  Such an equation naturally induces a
section $s$ of the bundle $\C\oplus S^4Q\oplus S^8Q$.  Consider a point $C\in
M^o$.
We abuse notation by allowing $C$ to also denote the corresponding curve.
Let $(q_0,q_1)$ be homogeneous coordinates on $\P^1$.  Identifying $C$ with
$\P^1$, we may describe
$C$ by equations of the form
\begin{equation}
\label{param}
x_0=f_4(q_0,q_1),\
x_i=l_i(q_0,q_1).
\end{equation}
 The equation $x_0-f_4(q_0,q_1)=0$
and its multiples for varying $C$ form the tautological
subbundle $\O_{\P}(-1)\subset\C\oplus S^4Q$ on $\P=\P(\C\oplus S^4Q^*)$.
$C$ is contained in $X$ if and only if an equation for $X$, when pulled back
to $\P^1$ via a parametrization of $C$, vanishes.  Substituting from the
second of equations~(\ref{param}) into (\ref{hyper}), it is seen that this
happens
if and only if $ax_0^2+g_4(l_1(q_0,q_1),\ldots,l_4(q_0,q_1))x_0+
g_8(l_1(q_0,q_1),\ldots,l_4(q_0,q_1))$ is a multiple of $x_0-f_4(q_0,q_1)$.
Multiplication induces an inclusion of bundles
$$(\C\oplus S^4Q)\otimes
\O_{\P}(-1)\hookrightarrow \C\oplus S^4Q\oplus S^8Q.$$
Putting all this together, we
see that $C\subset X$ if and only the section ${\bar s}$
of
$$B=(\C\oplus S^4Q\oplus S^8Q)
/((\C\oplus S^4Q)\otimes \O_{\P}(-1))$$
induced by
$s$ vanishes at $C$.  Note that if $C\in M-M^o$, then $C$ corresponds to a
curve
defined by equations of the form $f_4(q_0,q_1)=0,\ x_i=l_i(q_0,q_1)$.
Since such a curve would contain $p$, it follows that $C$ is not in the zero
locus of ${\bar s}$.  Also note that
$\dim(M)=rank(B)=9$; so one expects finitely many zeros of such a section;
hence finitely many weighted $\P^1$'s.  It is easy to prove that this is
indeed the case for general $X$.
The actual number is the degree of $c_9(B)$.  This may be calculated by
standard techniques in intersection theory \cite{fintthy} and the calculation
may be implemented via
\Schubert\ \cite{schub}.

The case $k=2$ is similar. Changing the meaning of the notation in the
obvious manner, one must consider $M=\P(\C\oplus S^2Q^*)$, and
calculate the degree of $c_7(B')$, where
$$B'=(\C\oplus S^2Q\oplus S^4Q\oplus S^6Q)
/((\C\oplus S^2Q\oplus S^4Q)\otimes\O_{\P}(-1)).$$

Combining these with the well-known number of lines on a quintic threefold,
the calculation of some examples considered in \cite{fontper,kt,morp-f} via
mirror symmetry may be verified.  The results are displayed in
table~\ref{table1}.

\begin{table}\begin{center}\begin{tabular}{|ccc|} \hline
Ambient space & Weighted degree & Number of lines \\ \hline
$\P(1^5)$ &  5 & 2875 \\
$\P(2,1^4)$ &  6 & 7884 \\
$\P(4,1^4)$ &  8 & 29504 \\ \hline
\end{tabular}\end{center}\caption{The number of lines.}
\label{table1}\end{table}

\bigskip\noindent
{\bf Problem:} Verify the predictions of mirror symmetry for weighted $\P^1$'s
in a weight 10 hypersurface in $\P(5,2,1^3)$.  Also, verify the predictions
of mirror symmetry for weighted conics on the weighted hypersurfaces
considered in this section.

\bigskip\noindent
{\em Remark:} The family of weighted conics on the general weighted octic in
$\P(4,1^4)$ is positive dimensional (independently observed by Koll{\'a}r);
hence part of the problem in this case
is to systematically assign numbers to positive dimensional families.  This
can be defined as the number of such curves that remain almost holomorphic
under a general almost complex deformation; but it is desirable to
give a purely algebraic description.

\section{Lines on higher dimensional varieties}
In this section and the next, we consider rational curves on
the generic Calabi-Yau hypersurface $X$ in $P^{k+1}$.  This is a
hypersurface of dimension $k$ and degree $k+2$.  For $k>3$, there will be
infinitely many lines and conics contained in $X$.  But there will only be
finitely many lines or conics which satisfy certain incidence properties
with fixed linear subspaces.

Since the normal bundle $N$ of $C$ in $P^{k+1}$ has degree $-2$, one expects
that for general $X$ and any $C\subset X$, $N\simeq\O\oplus\ldots\oplus\O
\oplus\O(-1)\oplus\O(-1)$ (with $k-3\ \O$'s).  Since $h^0(N)=k-3$ and
$h^1(N)=0$ in this case, the scheme of rational
curves on $X$ is expected to have dimension $k-3$.

For each $i$, let $L_i\subset\P^{k+1}$ denote a general linear subspace of
codimension $i$.  Pick positive integers $a,\ b,\ c$ such that $a+b+c=k$.
Following \cite{wittsm,morht,gmp}, define an invariant
$n^a_b(d)$ of $X$ as the number of holomorphic immersions
$f:\P^1\to X$ with $f(\P^1)$ of degree $d$
such that $f(0)\in L_a,\ f(1)\in L_b,\linebreak f(\infty)\in  L_c$.
These numbers, called ``Gromov-Witten invariants'' in \cite{morht},
are expected to be finite.  Note that the value of $c$ is implicit in the
notation $n^a_b(d)$ by virtue of the equation $a+b+c=k$.

These invariants are essentially the same as the number of reduced, irreducible
rational curves of degree $d$ in $X$ which meet each of $L_a,\ L_b$, and
$L_c$.  The $n^a_b(d)$ differ from the corresponding numbers of curves by one
factor of $d$ for each of the indices $a,b$, or $c$ equal to 1 (since
$C$ meets a general $L_1\ d$ times).  There is no difference for lines;
and for conics, we will see that in the calculation of the number of conics
satisfying the required incidence properties, the Gromov-Witten invariants
arise naturally.  So the Gromov-Witten invariants will be calculated and
tabulated, while the numbers of rational curves follow immediately
by division by the appropriate power of $d$, if necessary.

In the remainder of this section, we specialize to $d=1$, i.e. lines.
A theorem of Barth-van de Ven \cite{bv} states that
the Fano variety of lines on a degree $l$ hypersurface $X \subset \P^n$ is
smooth of dimension $2n-l-3$ for generic $X$
when $l+3 \le 2n$.
Applied in the present context of $X_{k+2}\subset\P^{k+1}$, we find
that the variety of lines must be smooth of dimension $k-3$ whenever
$k\ge3$.
{}From this, standard techniques show that a general $X$ contains finitely many
lines which meet each of $L_a,\ L_b$, and $L_c$.

So we can calculate the Gromov-Witten invariants by using the
Schubert calculus.  The lines are parametrized by the Grassmannian
$G(2,k+2)$.  The class of lines meeting $L_a$ is the Schubert cycle
\s{a-1}; similarly for $L_b$ and $L_c$.
Let $Q$ be the rank 2 universal quotient bundle on $G$.
Since the class of the variety of lines on $X$ is represented by
$c_{k+3}(S^{k+2}Q)$ and dimensions work out correctly,
the answer is the degree of
$c_{k+3}S^{k+2}Q\cdot\s{a-1}\cdot\s{b-1}\cdot\s{c-1}$.
These may be easily worked out as integers using \Schubert.
The answers obtained are displayed in table~\ref{table3}.

\bigskip
The original predictions for the numbers found in \cite{gmp}
resulted from a two-step process arising from mirror symmetry and conformal
field theory.  First, the $n^1_b(d)$ are found, followed by what amounts to
an expression for any $n^a_b(d)$ in terms of the various $n^1_{b'}(d')$ for
$d'\le d$.  Most of these expressions remain a mathematical mystery at
present.
However, the case $d=1$ can be established mathematically as follows.

\begin{fact}
Let $X$ be any Calabi-Yau manifold of dimension $k$ in any projective space.
Define $n^a_b(1)$ as above.  Assume that there are
finitely many lines in $X$ satisfying each of the respective incidence
conditions needed to define the $n^a_b(1)$.  Then
$$n^i_j(1)=\sum_{l=0}^{j-1}n^1_{i+l}(1)-\sum_{l=1}^{j-1}n^1_l(1).$$
\end{fact}

\bigskip\noindent
{\em Proof (sketch).}  Follows immediately by intersecting the cycle class
(in the
appropriate Grassmannian) of the
scheme of lines in $X$ with the identity
$$\s{i-1}\s{j-1}\s{k-i-j-1}=\sum_{l=0}^{j-1}\s{i+l-1}\s{k-i-l-2}-
\sum_{l=1}^{j-1}\s{l-1}\s{k-l-2},$$
an identity which can be proven by a few applications of Pieri's formula.

\begin{table}\small\begin{center}\begin{tabular}{|c|c|} \hline
  $k$ & $n^a_b(1)$\\ \hline
  3   & $n^1_1(1)=2875$ \\ \hline
  4   & $n^1_1(1)=60480$ \\ \hline
  5   & $n^1_1(1)=1009792,\ n^1_2(1)=1707797$ \\ \hline
  6   & $n^1_1(1)=15984640,\ n^1_2(1)=37502976,\ n^2_2(1)=59021312$ \\ \hline
  7   & $n^1_1(1)=253490796,\ n^1_2(1)=763954092,\ n^1_3(1)=1069047153$ \\
\cline{2-2}
  & $n^2_2(1)=1579510449$ \\ \hline
  8   & $n^1_1(1)=4120776000,\ n^1_2(1)=15274952000,\ n^1_3(1)=27768048000$ \\
\cline{2-2}
  & $n^2_2(1)=38922224000,\ n^2_3(1)=51415320000$ \\ \hline
  9   & $n^1_1(1)=69407571816,\ n^1_2(1)=307393401172,\ n^1_3(1)=695221679878$
\\ \cline{2-2}
  & $n^1_4(1)=905702054829,\ n^2_2(1)=933207509234,\ n^2_3(1)=1531516162891$ \\
\cline{2-2}
  & $n^3_3(1)=1919344441597$ \\ \hline
  10  & $n^1_1(1)=1217507106816,\ n^1_2(1)=6306655500288$ \\ \cline{2-2}
  & $n^1_3(1)=17225362851840,\ n^1_4(1)=28015971489792$ \\ \cline{2-2}
  & $n^2_2(1)=22314511245312,\ n^2_3(1)=44023827234816$ \\ \cline{2-2}
  & $n^2_4(1)=54814435872768,\ n^3_3(1)=65733143224320$ \\ \hline
\end{tabular}\end{center}\caption{Gromov-Witten invariants for lines.}
\label{table3}\end{table}

\section{Conics on higher dimensional varieties}
It can easily be shown that if $X$ is a general
hypersurface of degree $k+2$ in $\P^{k+1}$, then the variety of conics on $X$
has the expected dimension $k-3$.  Standard techniques show that
given positive integers $a,b,c$ with
$a+b+c=k$, there will be a finite number of conics in $X$ which meet each of
$L_a,\ L_b$, and $L_c$.  Thus the Gromov-Witten invariants are finite.  They
will be calculated here; the answers obtained are displayed in
table~\ref{table4}.


\begin{table}\footnotesize\begin{center}\begin{tabular}{|c|c|} \hline
  $k$& $n^a_b(2)$ \\ \hline
  3 & $n^1_1(2)=4874000$ \\ \hline
  4 & $n^1_1(2)=1763536320$ \\ \hline
  5 & $n^1_1(2)=488959144352,\ n^1_2(2)=1021575491286$ \\ \hline
  6 & $n^1_1(2)=133588638826496,\ n^1_2(2)=448681408315392 \
n^2_2(2)=821654025830400$ \\ \hline
  7 & $n^1_1(2)=39031273362637440,\ n^1_2(2)=187554590257349088$ \\ \cline{2-2}
  & $n^1_3(2)=312074852318965368,\ n^2_2(2)=506855012110118424$ \\ \hline
  8 & $n^1_1(2)=12607965435718224000,\ n^1_2(2)=80684596772238448000$ \\
\cline{2-2}
  & $n^1_3(2)=200581960800610752000,\ n^2_2(2)=295035175517918176000$ \\
\cline{2-2}
  & $n^2_3(2)=444475303469701680000$\\ \hline
  9 & $n^1_1(2)=4565325719860021608624,\ n^1_2(2)=37005001823802188657624$ \\
\cline{2-2}
  & $n^1_3(2)=127922335050535174614916,\ n^1_4(2)=193693669320390878077186$ \\
\cline{2-2}
  & $n^2_2(2)=173901546566279203106468,\ n^2_3(2)=364629304647788940660824$ \\
\cline{2-2}
  & $n^3_3(2)=498705676383823268404990$
  \\ \hline
  10 & $n^1_1(2)=1861791822397620935737344,\
n^1_2(2)=18415607624138339954786304$ \\ \cline{2-2}
  & $n^1_3(2)=83885220561474498867757056,\
n^1_4(2)=179982840924749584358866944$ \\ \cline{2-2}
  & $n^2_2(2)=107227899142191919158312960,\
n^2_3(2)=297755098999730079369412608$ \\ \cline{2-2}
  & $n^2_4(2)=417950364467570984815214592,\
n^3_3(2)=527556832251612742800359424$ \\ \hline
\end{tabular}\end{center}\caption{Gromov-Witten invariants for conics.}
\label{table4}\end{table}

We start with the well known description of the moduli space of conics in
$\P^{k+1}=\P(V)$, where $V$ is a $k+2$-dimensional vector space.  To describe
a conic, we first describe the 2 plane it spans,
and then choose a quadric in that 2-plane (up to scalar).  So let
$G=G(3,V)$ be the Grassmannian of 2-planes in $\P(V)$ (that is, of rank 3
quotients of $V$), and let $Q$ be
the universal rank 3 quotient bundle of linear forms on the varying subspace.
Then the moduli space of conics is $M=\P(S^2Q^*)$.  Following the reasoning in
section~\ref{weight-section} (or \cite{finite}), the scheme of conics on $X$ is
given by the locus over which a certain section of $F=S^{k+2}Q/(S^kQ
\otimes\O_{\P}(-1))$ vanishes.  Here $\O_{\P}(1)$ is the tautological sheaf
on $\P(S^2Q^*)$.  Since $F$ has rank $2k+5$, the conics on $X$ are represented
by $c_{2k+5}(F)$.

It remains to find the condition that a conic $C$ meets $L_a$.
One way to find this is to consider the moduli
space ${\cal M}$ of pointed conics, i.e. pairs $(p,C)$, with $C$ a conic, $p\in
C$.
This may easily be constructed as a bundle over $\P^{k+1}$, with fiber
over $p\in \P^{k+1}$ being the set of conics containing $p$.  We start by
constructing the moduli space of pointed 2-planes as follows.  Consider the
tautological exact sequence on $\P^{k+1}$:

$$0\to K\to V_{\P^{k+1}}\to \O(1)\to 0,$$

\noindent
where $V_{\P^{k+1}}$ is a trivial bundle of rank $k+2$ on $\P^{k+1}$
(more generally, $E_Y$
will stand for the pullback of $E$ to $Y$, the morphism used for the pullback
assumed to be clear in context).  Let $H=G(2,K)$ be the
Grassmannian of rank 2 quotients of $K$, \Q\ its universal rank 2
quotient, and ${\cal S}\subset K_H$ the universal subbundle.  These
fit into the exact sequence
$$0\to{\cal S}\to K_H\to \Q\to 0$$
of sheaves on $H$.  The natural
quotient $V_H\to V_H/{\cal S}$ induces a map $H\to G$ by the universal property
of the Grassmannian; since $V_H\to \O(1)_H$ clearly factors through
$V_H/{\cal S}$, it is easy to see that $H$ may be identified with the space of
pointed 2-planes in $\P(V)$.  Here $\O(1)$ denotes the tautological sheaf on
$\P(V)$ as before.

The conics containing $p$ globalize to a rank 5 bundle $W$ on $H$.  This
bundle is in fact the the kernel of the natural map
$S^2(V_H/{\cal S})\to\O(2)_H$.
Then the moduli space $M'$ of pointed conics may be seen to be
$\P(W^*)$.  Let $\O_W(1)$ be its tautological bundle.

Let $h=c_1(\O(1)_{M'})$.  Consider the natural morphism
$f:M'\to M$.  The variety of conics meeting $L^a$ is represented by the
class $f_*(h^a)$ for $a>1$.  Note that for $a=1$, $f_*(h)=2$.  This factor
exactly gives the factor needed to give the Gromov-Witten invariants rather
than the number of conics meeting three linear subspaces.  So the Gromov-Witten
invariants are given by the formula
$n^a_b(2)=\int_Mc_{2k+5}(Q)f_*(h^a)f_*(h^b)f_*(h^c)$, which is valid since
the dimensions work out correctly.

To compute these as numbers using \Schubert, everything is clear, except the
description of the morphism $f$.  But this may be described merely by knowing
the pullbacks $f^*(Q)$ and $f^*(\O_{\P}(1))$.  However, from the above
description and the universal properties, this is just $V_H/{\cal S}$ and
$\O_W(1)$.  \Schubert\ takes care of the rest.




\begin{thebibliography}{10}

\bibitem{bv}
W.~Barth and A.~van~de Ven.
\newblock {F}ano-varieties of lines on hypersurfaces.
\newblock {\em Arch.\ Math.}, 31(1):96--104, 1978.

\bibitem{cogp}
P.~Candelas, X.C. de~la Ossa, P.S. Green, and L.~Parkes.
\newblock A pair of {C}alabi-{Y}au manifolds as an exactly soluble
  superconformal theory.
\newblock {\em Nucl. Phys. B}, 359:21--74, 1991.

\bibitem{escub}
G.~Ellingsrud and S.~A. Str{\o}mme.
\newblock The number of twisted cubic curves on the general quintic threefold.
\newblock Univ.\ of Bergen Preprint 63-7-2-1992.

\bibitem{fontper}
A.~Font.
\newblock Periods and duality symmetries in {C}alabi-{Y}au compactifications.
\newblock Universidad Central de Venezuela preprint UCVFC/DF-1-92, 1992.

\bibitem{fintthy}
William Fulton.
\newblock {\em Intersection Theory}.
\newblock Springer-Verlag, Berlin Heidelberg New York Tokyo, 1984.

\bibitem{gmp}
B.~R. Greene, D.~R. Morrison, and M.~R. Plesser.
\newblock Mirror manifolds in higher dimension.
\newblock In preparation.

\bibitem{finite}
Sheldon Katz.
\newblock On the finiteness of rational curves on quintic threefolds.
\newblock {\em Comp.\ Math.}, 60:151--162, 1986.

\bibitem{schub}
Sheldon Katz and Stein~Arild Str{\o}mme.
\newblock {\sc schubert}: a {\sc maple} package for intersection theory.
\newblock Available by anonymous ftp from ftp.math.okstate.edu or
  linus.mi.uib.no, cd pub/schubert.

\bibitem{kt}
A.~Klemm and S.~Thiesen.
\newblock Considerations of one-modulus {C}alabi-{Y}au compactifications:
  {P}icard-{F}uchs equations, {K}{\"a}hler potentials and mirror maps.
\newblock Universit{\"a}t Karlsruhe and Technische Universit{\"a}t M{\"u}nchen
  preprint KA-THEP-03/92, TUM-TH-143-92, 1992.

\bibitem{ltcy}
A.~Libgober and J.~Teitelbaum.
\newblock Lines on {C}alabi {Y}au complete intersections, mirror symmetry, and
  {P}icard {F}uchs equations.
\newblock Preprint, University of Illinois at Chicago.

\bibitem{morp-f}
D.~R. Morrison.
\newblock Picard-{F}uchs equations and mirror maps for hypersurfaces.
\newblock In S.-T. Yau, editor, {\em Essays on Mirror Manifolds}, pages
  241--264. International Press, Hong Kong, 1992.

\bibitem{morht}
D.~R. Morrison.
\newblock {H}odge-theoretic aspects of mirror symmetry.
\newblock In preparation.

\bibitem{wittsm}
E.~Witten.
\newblock Topological sigma models.
\newblock {\em Commun.\ Math.\ Phys.}, 118:411--449, 1988.

\bibitem{yau}
S.-T. Yau, editor.
\newblock {\em Essays on Mirror Manifolds}, Hong Kong, 1992. International
  Press.

\end{thebibliography}
\end{document}